\documentclass{article}
\usepackage{amsfonts}
\usepackage{amssymb}
\usepackage{graphicx}
\usepackage{amsmath}

\newtheorem{theorem}{Theorem}
\newtheorem{acknowledgement}[theorem]{Acknowledgement}

\newtheorem{condition}[theorem]{Condition}

\newtheorem{corollary}[theorem]{Corollary}

\newtheorem{definition}[theorem]{Definition}

\newtheorem{lemma}[theorem]{Lemma}

\newtheorem{proposition}[theorem]{Proposition}
\newtheorem{remark}[theorem]{Remark}

\newenvironment{proof}[1][Proof]{\textbf{#1.} }{\ \rule{0.5em}{0.5em}}
\input{tcilatex}

\begin{document}

\title{A Note on the Notion of Geometric Rough Paths}
\author{Peter Friz\thanks{%
Courant Institute, 251 Mercer St. New York, NY 10012, USA,
Peter.Friz@cims.nyu.edu} and Nicolas Victoir$\thanks{%
Magdalen College, Oxford OX1 4AU, UK; victoir@maths.ox.ac.uk}$}
\maketitle

\begin{abstract}
We use simple sub-Riemannian techniques to prove that an arbitrary geometric 
$p$-rough path in the sense of \cite{Ly} is the limit in sup-norm of a
sequence of canonically lifted smooth paths, which are uniformly bounded in $%
p$-variation, clarifying the two different defintions of a geometric $p$%
-rough \cite{Ly,LQ}.

Our proofs are based on fine estimates in terms of control functions and are
sufficiently general to include the case of H\"{o}lder- and modulus-type
regularity \cite{F,FV}. This allows us to extend a few classical results on H%
\"{o}lder-spaces \cite{Ci,MS} and $p$-variation spaces \cite{Du,Wi} to the
non-commutative setting necessary for the theory of rough paths.
\end{abstract}

\section*{Introduction}

Over the last years T. Lyons developed a general theory of integration and
differential equations of the form 
\begin{equation}
dy_{t}=f(y_{t})dx_{t}  \label{ODE}
\end{equation}%
where the driving signal $x_{t}\in V,$ $t\in \lbrack 0,1],$ is a Banach
space valued path of finite $p$-variation. For $p<2$ this leads to a
(pathwise) differential equation theory based on Young integrals, but it was
not observed before \cite{L94} that $x_{.}\mapsto y_{.}$ is actually
continuous in $p$-variation topology (and also in some more refined
topologies). Most recently \cite{LL}, Fr\'{e}chet smoothness was established
under natural conditions on $f\,$. The situation is much more complicated
for $p\geq 2$ and was successfully worked out in \cite{Ly}. Clearly, this
case is the one needed for applications to stochastic differential
equations. The path $x$ driving the differential equation (\ref{ODE}) needs
to be lifted to a path $X$ of finite $p$-variation with values in $%
G^{[p]}(V) $, the free nilpotent group of step $[p]$ over $V$. The theory of
rough paths then gives a solution $y$ to the differential equation, and
actually also automatically lifts $y$ to a path $Y$ of finite $p$-variation
with values in a free nilpotent group of step $[p]$. Moreover, the map $%
X\rightarrow Y$ is continuous using an appropriate $p$variation distance
(and actually also in more refined topologies).

A smooth $V$-valued path $x$ can be canonically lifted to a $G^{[p]}(V)$%
-valued path $X$. The solution of equation (\ref{ODE}) driven by such a
canonical $X$ is simply the canonical lift of the classical solution $y$ of
the corresponding ODE. If $x$ is a Brownian motion and $X$ its Stratonovich
lift to a geometric $p$-rough path, then $y$ is the solution of the
corresponding Stratonovich SDE.

The theory of rough paths tells us that the signal in control differential
equations of type (\ref{ODE}) are paths with values in a free nilpotent
group, satisfying some $p$-variation constraints. The set of such signals is
called the set of geometric $p$-rough paths.

There has been some confusions on the precise definition of a geometric $p$%
-rough path: in \cite{Ly}, a geometric $p$-rough path is defined as the set
of paths with values in $G^{[p]}(V)$ which has finite $p$-variation,
computed with a natural metric associated to the group. It is then falsely
claimed that an equivalent definition is the $p$-variation closure of the
canonical lift of smooth paths to paths with values in the group. The latter
definition was the one chosen in the more recent monograph of Lyons and Qian 
\cite{LQ} and we shall follow its notation.

\bigskip

This paper studies precisely the difference between these two definitions.
In the first section, we reintroduce the algebra and analysis needed to
explain the theory of rough paths:\ free nilpotent groups and their
homogeneous norms. The second section deals with some basic results on path
space. There are at least two possible notions of generalization of H\"{o}%
lder distance between two group valued paths. We show that these two notions
lead to the same topology. Then, we obtain some classical interpolation
results. The third and fourth section study precisely the set of paths with
values in $G^{[p]}(V)$ which have finite $p$-variation, and the set of
geometric $p$-rough paths. We will see that if $Y$ is a $G^{[p]}(V)$-valued
path with finite $p$-variation, it is the limit in sup-norm of a sequence
uniformly bounded in $p$-variation norm of signature of smooth paths. These
smooth paths are constructed using (almost) sub-riemannian geodesics \cite%
{Mo}. If $Y$ is indeed a geometric $p$-rough path, this sequence is shown to
converge in $p$-variation distance. The same results holds replacing $p$%
-variation by $1/p$-H\"{o}lder, or some more general modulus norms.

We also give a characterization of the set of geometric $p$-rough paths in
the spirit of the Wiener class \cite{Du,Wi}. When translated into $1/p$-H%
\"{o}lder topology, the Wiener class relates to a characterization due to
Ciesielski in the vector space case \cite{Ci,MS}. Finally, we precise which
of the spaces under consideration are Polish.

$C$ in this paper denotes a constant, which may vary from line to line.

\section{Algebraic Preliminaries}

We refer to \cite{Re} for more details on free nilpotent groups, and \cite%
{FS,Mo} on homogeneous norms and Carnot Caratheodory distance.

\subsection{Free Nilpotent Groups}

We fix a real Banach space $\left( V,\left\Vert .\right\Vert \right) $, that
we assume finite dimensional. Let $T(V)=\bigoplus_{n=0}^{\infty }V^{\otimes
n}$ be the tensor algebra over $V$. $T(V)$ equipped with standard addition $%
+ $, tensor multiplication $\otimes $ and scalar product is an associative
algebra. $T^{(m)}(V),$ the quotient algebra of $T(V)$ by the ideal $%
\bigoplus_{n=m+1}^{\infty }V^{\otimes n},$ inherits this algebra structure.
\bigskip One can define on $T^{(m)}(V)$ a Lie bracket by the formula 
\begin{equation*}
\lbrack a,b]=a\otimes b-b\otimes a,
\end{equation*}%
which makes $T^{(m)}(V)$ into a Lie algebra. Let $\mathcal{G}^{m}(V)$ be the
Lie subalgebra of $T^{(m)}(V)$ generated by elements in $V$. Note that 
\begin{equation*}
\mathcal{G}^{m}(V)\simeq \bigoplus_{i=1}^{m}V_{i},
\end{equation*}%
where 
\begin{equation}
V_{1}=V\text{ and }V_{i+1}=[V,V_{i}].  \label{defV_i}
\end{equation}%
$\mathcal{G}^{m}(V)$ is the free nilpotent Lie algebra of step $m$ \cite%
{Ly,LQ,Re}. The exponential, logarithm and inverse function are defined on $%
T^{(m)}(V)$ by mean of their power series. We denote by $G^{m}(V)=\exp
\left( \mathcal{G}^{m}(V)\right) .$ By the Baker-Campbell-Hausdorff formula, 
$\left( G^{m}(V),\otimes \right) $ is a connected nilpotent Lie group,
called the free nilpotent Lie group of step $m$ over $V$, with Lie algebra $%
\mathcal{G}^{m}(V)$. We also define $\widetilde{T}^{(m)}(V)$ to be the set
of elements in $\widetilde{T}^{(m)}$ such that the term in $V^{\otimes 0}=%
\mathbb{R}$ is equal to $1$. $\widetilde{T}^{(m)}(V)$ with the product $%
\otimes $ of $T^{(m)}(V)$ is a Lie group. Note that $G^{m}(V)$ is a subgroup
of $\widetilde{T}^{(m)}(V)$. For an element $g=1+v_{1}+\ldots +v_{m}\in 
\widetilde{T}^{(m)}(V)\,$, with $v_{i}\in V^{\otimes i}$, we define, for $%
t\in \mathbb{R}$,%
\begin{equation*}
\delta _{t}g=1+tv_{1}+\ldots +t^{m}v_{m}.
\end{equation*}%
$\delta $ is called the dilation operator.

\subsection{Homogeneous Norms}

We are now going to equip $G^{m}(V)$ with a (symmetric sub-additive)
homogeneous norm \cite{FS}, i.e. a function $\left\Vert .\right\Vert
_{G^{m}(V)}:G^{m}(V)\rightarrow \mathbb{R}^{+}$ such that 
\begin{equation*}
\begin{tabular}{ll}
(i) & $\left\Vert .\right\Vert _{G^{m}(V)}$ if and only if $g=1$, \\ 
(ii) & $\left\Vert \delta _{t}g\right\Vert _{G^{m}(V)}=\left\vert
t\right\vert \left\Vert g\right\Vert _{G^{m}(V)},$ \\ 
(iii) & for all $g,h\in _{G^{m}(V)}$, $\left\Vert g\otimes h\right\Vert
_{G^{m}(V)}\leq \left\Vert g\right\Vert _{G^{m}(V)}+\left\Vert h\right\Vert
_{G^{m}(V)},$ \\ 
(iv) & for all $g$, $\left\Vert g\right\Vert _{G^{m}(V)}=\left\Vert
g^{-1}\right\Vert _{G^{m}(V)}.$%
\end{tabular}%
\newline
\end{equation*}%
We define on the group the Carnot-Caratheodory homogeneous norm $%
||.||_{G^{m}(V)}$ with the help of the formula%
\begin{equation*}
||g||_{G^{m}(V)}=\inf \left( \int_{0}^{1}\left\vert \dot{y}_{r}\right\vert
dr\right) ,
\end{equation*}%
where the infimum is taken over all smooth paths $y:\left[ 0,1\right]
\rightarrow V$ such that%
\begin{equation}
S_{m}(y)_{0,1}=g.  \label{chow}
\end{equation}%
Here $S_{m}$ denotes the $m$-signature of $y$ in $G^{m}(V)$ between the time 
$r$ and $s$, that is%
\begin{equation*}
S_{m}(y)_{r,s}=\left(
1,y_{r,s}=\int_{r}^{s}dy_{u},\int_{r}^{s}y_{r,u}\otimes
dy_{u},...,\int_{r<u_{1}<\cdots <u_{m}<s}dy_{u_{1}}\otimes ...\otimes
dy_{u_{m}}\right) .
\end{equation*}%
The fact that there exists a smooth path $y$ which satisfies (\ref{chow}) is
precisely Chow's theorem \cite{Mo}. Chen's theorem \cite{Ch} asserts that $%
S_{m}(y)_{0,r}\otimes S_{m}(y)_{r,s}=S_{m}(y)_{0,s}$. Note that $%
s\rightarrow S_{m}(y)_{0,s}$ is equivalently defined as the solution of the
ordinary differential equation in $T^{(m)}(V)$%
\begin{equation*}
dS_{m}(y)_{0,s}=S_{m}(y)_{0,s}\otimes dy_{s}.
\end{equation*}

\begin{proposition}
\label{obvProp}Let $z$ be a path $[0,1]\rightarrow V$ in $W^{1,1}$ (i.e.
with derivatives in $L^{1}$). Then%
\begin{equation*}
\left\Vert S_{m}(z)_{s,t}\right\Vert _{G^{m}(V)}\leq \int_{s}^{t}\left\vert 
\dot{z}\right\vert dr.
\end{equation*}
\end{proposition}

\begin{proof}
Obvious by the definition of the Carnot-Caratheodory norm.
\end{proof}

We now fix $\left\vert .\right\vert _{i}$ be some norms on $V^{\otimes i}$
such that for all $\left( a^{i},a^{j}\right) \in V^{\otimes i}\times
V^{\otimes j}$, $\left\vert a^{i}\otimes a^{j}\right\vert _{i+j}\leq
\left\vert a^{i}\right\vert _{i}+\left\vert a^{j}\right\vert _{j}$. To
simplify notations, we will write $\left\vert .\right\vert $ for all these
norms. For $x\in \widetilde{T}^{(m)}(V)$, 
\begin{equation*}
\left\Vert x\right\Vert _{\widetilde{T}^{(m)}(V)}=\max_{i=1,\ldots ,m}\left(
i!\left\vert x^{i}\right\vert \right) ^{1/i},
\end{equation*}%
where $x=1+x^{1}+\ldots +x^{m},$ $x^{i}\in V^{\otimes i}.$ Then $x\in 
\widetilde{T}^{(m)}(V)\rightarrow \left\Vert x\right\Vert _{\widetilde{T}%
^{(m)}(V)}$ defines a subadditive homogeneous norm on $\widetilde{T}%
^{(m)}(V) $ \footnote{%
Note that $g\in \widetilde{T}^{(m)}(V)\rightarrow \left\Vert g\right\Vert _{%
\widetilde{T}^{(m)}(V)}+\left\Vert g^{-1}\right\Vert _{\widetilde{T}%
^{(m)}(V)}$ defines a subadditive symmetric homogeneous norm, which is
equivalent to $\left\Vert .\right\Vert _{\widetilde{T}^{(m)}(V)}$. Indeed,\
if $g\in \widetilde{T}^{\left( m\right) }(V),$ $g=g^{0}+g^{1}+...+g^{m}$, $%
g^{i}\in V^{\otimes i}$, (with $g^{0}=1$), then, for $k\geq 1$%
\begin{equation}
\left( g^{-1}\right) ^{k}=\sum_{j=1}^{k}(-1)^{j}\sum_{\substack{ %
i_{1},\cdots ,i_{j}\in \{1,\cdots ,m\}  \\ i_{1}+\cdots +i_{j}=k}}%
g_{i_{1}}\otimes \cdots \otimes g_{i_{j}}.  \label{inverse}
\end{equation}%
This easily implies that there exists a constant $C_{m}$,which depends only
on $m$, such that%
\begin{equation*}
\left\Vert g^{-1}\right\Vert _{\widetilde{T}^{(m)}(V)}\leq C_{m}\left\Vert
g\right\Vert _{\widetilde{T}^{(m)}(V)}\text{.}
\end{equation*}%
}. When restricted to $G^{m}(V)$, it is also symmetric \cite{Ly}. For $%
g=\exp (\ell ^{1}+\ldots +\ell ^{m})$, with $\ell _{i}\in \mathcal{L}_{i}$,
we also define%
\begin{equation*}
\left\Vert g\right\Vert _{\mathcal{L}^{(m)}(V)}=\max_{i=1,\ldots
,m}\left\vert \ell ^{i}\right\vert ^{1/i}.
\end{equation*}

$\left\Vert .\right\Vert _{\mathcal{L}^{(m)}(V)}$ is a symmetric homogeneous
norm which is equivalent, even when $V$ is of infinite dimension, to the
homogeneous norm $\left\Vert .\right\Vert _{\widetilde{T}^{(m)}(V)}$
restricted to the group $G^{m}(V)$ \cite{LV}. When $V$ is finite
dimensional, as all homogeneous norms are equivalent \cite{Go}, $\left\Vert
.\right\Vert _{G^{m}(V)}$, $\left\Vert .\right\Vert _{\widetilde{T}%
^{(m)}(V)} $ restricted to the group $G^{m}(V)$, and $\left\Vert
.\right\Vert _{\mathcal{L}^{(m)}(V)}$ are equivalent. Therefore, when no
confusion arises, we will not distinguish between these homogeneous norms
and we will denote them $\left\Vert .\right\Vert _{m}$. We define a left
invariant distance: $d_{m}(g,h)=\left\Vert g^{-1}\otimes h\right\Vert _{m}.$

\section{Group Valued Paths}

By a $G$-valued path, where $\left( G,\cdot \right) $ is a Lie group, we
will always mean a continuous function from $[0,1]$ into $G$, starting at
the neutral element of the group. We denote this set by $C_{0}\left(
[0,1],G\right) $. Moreover, if $Y$ is such a path, we will use the \emph{%
notation} throughout the paper $Y_{s,t}=Y_{s}^{-1}\cdot Y_{t}$.

\subsection{Some Metrics giving the Same Topology}

\begin{lemma}
Let $g,h$ be two elements of $\widetilde{T}^{\left( m\right) }(V),$ with $%
g=1+g^{1}+...+g^{m}$, $g^{i}\in V^{\otimes i}$; we use similar notations for 
$h.$ The following equation holds in $V^{\otimes k},k=1,...,m$%
\begin{equation}
\left( g^{-1}\otimes h\right) ^{k}=h^{k}-g^{k}+\sum_{i=1}^{k-1}\left(
g^{-1}\right) ^{k-i}\otimes (h^{i}-g^{i})  \label{linkminusinverse}
\end{equation}
\end{lemma}

\begin{proof}
Set $g^{0}=h^{0}=1.$ By definition of the tensor product in $\widetilde{T}%
^{\left( m\right) }(V),$ $\left( g^{-1}\otimes h\right)
^{k}=\sum_{i=0}^{k}\left( g^{-1}\right) ^{k-i}\otimes h^{i}.$The result
follows from subtracting to the previous expression $0=\left( g^{-1}\otimes
g\right) ^{k}=\sum_{i=0}^{k}\left( g^{-1}\right) ^{k-i}\otimes g^{i}.$
\end{proof}

\begin{proposition}
\label{topo}Let $\varepsilon \in (0,1).$ Given $g\in \widetilde{T}^{(m)}(V)$
there exists a constant $C_{m}>0$ such that:\newline
(i) If $\max_{i=1,...,m}\left\vert h^{i}-g^{i}\right\vert \leq 1$,%
\begin{equation}
d_{m}(g,h)\leq C_{m}\max \left\{ 1,\left\Vert g\right\Vert _{m}\right\}
\left( \max_{i=1,...,m}\left\vert h^{i}-g^{i}\right\vert \right) ^{1/m}.
\label{firstine}
\end{equation}%
(ii)If $\left\Vert g^{-1}\otimes h\right\Vert _{m}\leq 1$,%
\begin{equation}
\max_{i=1,...,m}\left\vert h^{i}-g^{i}\right\vert \leq C_{m}\max \left\{
1,\left\Vert g\right\Vert _{m}^{m}\right\} d_{m}(g,h).  \label{secondine}
\end{equation}
\end{proposition}

\begin{proof}
From formula (\ref{linkminusinverse}), one easily sees that%
\begin{equation*}
\left\vert \left( g^{-1}\otimes h\right) ^{k}\right\vert \leq m\max \left\{
1,\left\Vert g^{-1}\right\Vert _{m}^{k}\right\} \max_{i=1,...,m}\left\vert
h^{i}-g^{i}\right\vert .
\end{equation*}%
Hence, 
\begin{eqnarray*}
\left\vert \left( g^{-1}\otimes h\right) ^{k}\right\vert ^{1/k} &\leq
&C_{m}\max \left\{ 1,\left\Vert g\right\Vert _{m}\right\} \left(
\max_{i=1,...,m}\left\vert h^{i}-g^{i}\right\vert \right) ^{1/k} \\
&\leq &C_{m}\max \left\{ 1,\left\Vert g\right\Vert _{m}\right\} \left(
\max_{i=1,...,m}\left\vert h^{i}-g^{i}\right\vert \right) ^{1/m}
\end{eqnarray*}%
which gives inequality (\ref{firstine}).\newline
Reciprocally, assume that $\left\Vert g^{-1}\otimes h\right\Vert _{m}\leq 1$%
. We are going to show by induction that there exists a constant $C_{m}$
such that for all $i\in \{1,\cdots ,m\}$, 
\begin{equation}
\left\vert h^{i}-g^{i}\right\vert \leq C_{m}\max \left\{ 1,\left\Vert
g\right\Vert _{m}^{i}\right\} \left\Vert g^{-1}\otimes h\right\Vert _{m}.
\label{whattoprove}
\end{equation}%
$h^{1}-g^{1}=\left( g^{-1}\otimes h\right) ^{1}$ so the initial step is
easy. Assume now that (\ref{whattoprove}) is true up to a fixed index $i$.
Inequality (\ref{linkminusinverse}) then gives%
\begin{eqnarray*}
\left\vert h^{i+1}-g^{i+1}\right\vert &\leq &\left\Vert g^{-1}\otimes
h\right\Vert _{m}^{k}+C_{m}\sum_{j=1}^{i}\left\Vert g\right\Vert
_{m}^{i+1-j}\left\vert h^{j}-g^{j}\right\vert \\
&\leq &\left\Vert g^{-1}\otimes h\right\Vert
_{m}+C_{m}\sum_{j=1}^{i}\left\Vert g\right\Vert _{m}^{i+1-j}\max \left\{
1,\left\Vert g\right\Vert _{m}^{j}\right\} \left\Vert g^{-1}\otimes
h\right\Vert _{m} \\
&\leq &C_{m}\left\Vert g^{-1}\otimes h\right\Vert _{m}\max \left\{
1,\left\Vert g\right\Vert _{m}^{i+1}\right\} .
\end{eqnarray*}%
A straight-forward modification of the above proof yields the same result in
terms of the right invariant distance.
\end{proof}

We obtain the following:

\begin{proposition}
\label{rightleft}Define a right invariant distance $d_{r,m}(g,h)=\left\Vert
g\otimes h^{-1}\right\Vert _{m}$ based on a homogeneous norm $\left\Vert
.\right\Vert _{m}$. If $g_{n}$, $n\in \mathbb{N}$, and $g$ are elements in $%
\widetilde{T}^{(m)}(V)$ then the following is equivalent\newline
(i):\ $\lim_{n\rightarrow \infty }d_{r,m}(g_{n},g)=0.$\newline
(ii): $\lim_{n\rightarrow \infty }\max_{i=1,...,m}\left\vert
g_{n}^{i}-g^{i}\right\vert =0.$\newline
(iii):\ $\lim_{n\rightarrow \infty }d_{m}(g_{n},g)=0.$
\end{proposition}

\begin{corollary}
\ Let $X,Y$ be $G^{m}(V)$-valued paths with $p$-variation controlled by $%
\omega $. Let $d$ denote $d_{m}$ or $d_{r,m}.$ There exists a constant $%
c=c\left( \left\Vert X\right\Vert _{\infty },m\right) $ such that for $%
\varepsilon $ small enough, namely $0<\varepsilon <1/\left( \omega (0,1)\vee
1\right) $, we have 
\begin{gather*}
\max_{k=1,\cdots ,m}\frac{\left\vert X_{s,t}^{k}-Y_{s,t}^{k}\right\vert }{%
\omega (s,t)^{k/p}}\leq \varepsilon \Longrightarrow d\left(
X_{s,t},Y_{s,t}\right) \leq c\varepsilon ^{1/m}\omega ^{1/p}(s,t) \\
d\left( X_{s,t},Y_{s,t}\right) \leq \varepsilon \omega
^{1/p}(s,t)\Longrightarrow \max_{k=1,\cdots ,m}\frac{\left\vert
X_{s,t}^{k}-Y_{s,t}^{k}\right\vert }{\omega (s,t)^{k/p}}\leq c\varepsilon 
\text{.}
\end{gather*}
\end{corollary}

\begin{proof}
Define using the dilation operator $\delta $ on $T^{(m)},$%
\begin{equation*}
\tilde{X}_{s,t}=\delta _{\gamma }(X_{s,t})\text{ with }\gamma =1/\omega
(s,t).
\end{equation*}%
As $\delta $ commutes with $\otimes $ (and $^{-1}$) so that%
\begin{eqnarray*}
d_{m}\left( \tilde{X}_{s,t},\tilde{Y}_{s,t}\right) &=&\left\Vert \delta
_{\gamma }(X_{s,t}^{-1}\otimes Y_{s,t})\right\Vert _{m} \\
&=&\gamma d_{m}\left( X_{s,t},Y_{s,t}\right)
\end{eqnarray*}%
This reduction allows us to consider without loss of generality $\omega
(s,t)=1$ and the proposition follows from the results above.
\end{proof}

The last corollary implies that the topologies induced by the distances on $%
G^{m}(V)$-valued path space%
\begin{equation*}
\sup_{0\leq s<t\leq 1}\max_{k=1,\cdots ,m}\frac{\left\vert
X_{s,t}^{k}-Y_{s,t}^{k}\right\vert }{\omega (s,t)^{k/p}}
\end{equation*}%
and%
\begin{equation*}
\sup_{0\leq s<t\leq 1}\frac{d_{m}\left( X_{s,t},Y_{s,t}\right) }{\omega
(s,t)^{1/p}}
\end{equation*}%
are the same (here $\omega (s,t)$ is a control\footnote{$\omega $ is a
control if%
\begin{equation}
\begin{tabular}{ll}
(i) & $\omega :\left\{ (s,t),0\leq s\leq t\leq 1\right\} \rightarrow \mathbb{%
\ R}^{+}$ is continuous. \\ 
(ii) & $\omega $ is super-additive, i.e. $\forall $ $s<t<u$, $\omega
(s,t)+\omega (t,u)\leq \omega (t,u)$. \\ 
(iii) & $\omega (t,t)=0$ for all $t\in \lbrack 0,1]$%
\end{tabular}
\label{control}
\end{equation}%
} equal to $0$ only on the diagonal). The first distance is the one used by
Lyons and Lyons/Qian for the continuity results of integration and It\^{o}
map, the second one is the authors' favorite one. Any continuity result can
therefore be stated in either distances.

\subsection{$p$-Variation and Modulus Distances}

A path $x$ in a $G^{m}(V)$ is said to have finite $p$-variation if for all
subdivision $D=(0=t_{0}<\cdots <t_{n}=1)$ of $[0,1]$,%
\begin{equation*}
\sum_{i=0}^{n-1}\left\Vert x_{t_{i},t_{i+1}}\right\Vert _{m}^{p}<\infty 
\text{.}
\end{equation*}%
It can easily be seen to be equivalent the existence of a control function $%
\omega $ such that for all $s\leq t$, $\left\Vert x_{s,t}\right\Vert
_{m}^{p}\leq \omega (s,t)$. We define the following metric on the space of $%
G^{m}(V)$-valued paths:%
\begin{equation*}
d_{\omega ,p}(x,y)=\sup_{0\leq s<t\leq 1}\frac{d_{m}(x_{s,t},y_{s,t})}{%
\omega (s,t)^{1/p}}.
\end{equation*}%
Note than when $\omega (s,t)=t-s$, $d_{\omega ,p}$ is just the $1/p$-H\"{o}%
lder distance. We introduce a class of \textquotedblleft
nice\textquotedblright\ controls:

\begin{condition}
A control is said to satisfy the condition $(H_{p})$ if it is not identical
equal to 0 and if there exists $C$ such that for all $r<s<t<u$, $\frac{%
\omega (r,u)}{\left( u-r\right) ^{p}}\leq C\frac{\omega (s,t)}{\left(
t-s\right) ^{p}}.$ Note that it implies in particular that $\left(
t-s\right) ^{p}\leq C\omega (s,t)$.
\end{condition}

The control $\left( s,t\right) \rightarrow t-s$, as well as the controls
introduced in \cite{FV}, satisfy condition $(H_{p})$.

We will also look at the $p$-variation distance:%
\begin{equation*}
d_{p-var}(x,y)=\sup_{D=(0=t_{0}<\cdots <t_{n}=1)}\left(
\sum_{i=0}^{n-1}d_{m}(x_{t_{i},t_{i+1}},y_{t_{i},t_{i+1}})^{p}\right) ^{1/p}.
\end{equation*}%
We define $\left\Vert x\right\Vert _{\omega ,p}=d_{\omega ,p}(x,0)$ and $%
\left\Vert x\right\Vert _{p-var}=d_{p-var}(x,0)$.

\begin{definition}
We define the following path-spaces%
\begin{eqnarray*}
C^{p-var}(G^{m}(V)) &=&\left\{ x\in C_{0}([0,1],G^{m}(V))\text{ such that }%
\left\Vert x\right\Vert _{p-var}<\infty \right\} , \\
C^{\omega ,p}(G^{m}(V)) &=&\left\{ x\in C_{0}([0,1],G^{m}(V))\text{ such
that }\left\Vert x\right\Vert _{\omega ,p}<\infty \right\} .
\end{eqnarray*}%
$C^{0,p-var}(G^{m}(V))$ (resp. $C^{0,\omega ,p}(G^{m}(V))$) is defined as
the $d_{p-var}$-closure (resp. $d_{\omega ,p}$-closure) of the set $\left\{
S_{m}(x),\text{ }x\text{ smooth }V\text{-valued path}\right\} $.
\end{definition}

$C^{0,p-var}(G^{[p]}(V))$ is precisely the set of geometric $p$-rough paths,
according to the definition of \cite{LQ}, while $C^{p-var}(G^{[p]}(V))$ is
the set of geometric $p$-rough paths, according to \cite{Ly}. There has
indeed been some confusions in the seminal paper \cite{Ly} between the two
sets $C^{0,p-var}(G^{[p]}(V))$ and $C^{p-var}(G^{[p]}(V)).$ Here, we propose
to study and characterize these sets precisely. Studying their subset $%
C^{0,\omega ,p}(G^{m}(V))$ and $C^{\omega ,p}(G^{m}(V))$ is also of
interest, as the continuity results of the theory of rough paths can involve
the distance $d_{\omega ,p}$.

We will need some interpolation results.

\subsection{Interpolations}

\begin{proposition}
\label{AA} Let $Y(n)$ be a sequence of equi-continuous $G^{m}(V)$-valued
paths converging pointwise to a continuous path $Y.$ Then $Y\left( n\right) $
converges uniformly on $[0,1]$ to $Y,$ i.e.%
\begin{equation*}
\sup_{t}d_{m}(Y(n)_{t},Y_{t})\rightarrow 0.
\end{equation*}
\end{proposition}

\begin{proof}
Standard Arzela-Ascoli argument.
\end{proof}

\begin{proposition}
$\widetilde{d}_{\infty }(Y(n),Y)=\sup_{t}d_{m}(Y_{t}(n),Y_{t})\rightarrow
_{n\rightarrow \infty }0$ if and only if 
\begin{equation*}
d_{\infty }(Y(n),Y)=\sup_{s,t}d_{m}(Y(n)_{s,t},Y_{s,t})\rightarrow
_{_{n\rightarrow \infty }}0.
\end{equation*}
\end{proposition}

\begin{proof}
Clearly, if $d_{\infty }(Y(n),Y)$ goes to $0$ as $n\rightarrow \infty ,$
then so does $\widetilde{d}_{\infty }(Y(n),Y)$. Reciprocally,%
\begin{equation*}
d_{\infty }(Y(n)_{s,t},Y_{s,t})\leq
\sup_{s,t}d_{m}(Y(n)_{s,t},Y(n)_{s}^{-1}\otimes
Y_{t})+\sup_{s,t}d_{m}(Y(n)_{s}^{-1}\otimes Y_{t},Y_{s,t}).
\end{equation*}%
But $\sup_{s,t}d_{m}(Y(n)_{s,t},Y(n)_{s}^{-1}\otimes
Y_{t})=\sup_{t}d_{m}(Y(n)_{t},Y_{t})$ goes, by assumption, to $0$ when $%
n\rightarrow \infty $. Moreover, by corollary \ref{rightleft}, 
\begin{equation*}
\lim_{n\rightarrow \infty }\sup_{s,t}d_{m}(Y(n)_{s}^{-1}\otimes
Y_{t},Y_{s,t})=0
\end{equation*}%
if and only if 
\begin{equation*}
\lim_{n\rightarrow \infty }\sup_{s,t}d_{r,m}(Y(n)_{s}^{-1}\otimes
Y_{t},Y_{s,t})=0.
\end{equation*}%
But the latter is true as 
\begin{eqnarray*}
d_{r,m}(Y(n)_{s}^{-1}\otimes Y_{t},Y_{s,t})
&=&d_{r,m}(Y(n)_{s}^{-1},Y_{s}^{-1}) \\
&=&d_{m}(Y(n)_{s},Y_{s}).
\end{eqnarray*}
\end{proof}

\begin{remark}
$\tilde{d}_{\infty }$ and $d_{\infty }$ are not equivalent distances, but
induce the same topology. The following inequalities are classical, at least
for the H\"{o}lder norms \cite{Kr,St} and $p$-variation norms \cite{Le}.
\end{remark}

\begin{proposition}
Let $1\leq p<p^{\prime }<\infty $. Then for all $G^{m}(V)$-valued paths $Y,Z$%
\begin{equation}
d_{\omega ,p^{\prime }}(Y,Z)\leq d_{\infty }(Y,Z)^{1-p/p^{\prime }}d_{\omega
,p}(Y,Z)^{\frac{p}{p^{\prime }}}.  \label{interpolation}
\end{equation}%
In particular, if $Y(n)$ converges pointwise to $Y$ and $\sup_{n}\left\Vert
Y(n)\right\Vert _{\omega ,p}<\infty $ then 
\begin{equation*}
d_{\omega ,p^{\prime }}(Y(n),Y)\rightarrow 0.
\end{equation*}
\end{proposition}

\begin{proof}
For all $s<t$,%
\begin{eqnarray*}
\frac{d_{m}(Y_{s,t},Z_{s,t})^{p^{\prime }}}{\omega (s,t)}
&=&d_{m}(Y_{s,t},Z_{s,t})^{p^{\prime }-p}\frac{d_{m}(Y_{s,t},Z_{s,t})^{p}}{%
\omega (s,t)} \\
&\leq &d_{\infty }(Y,Z)^{p^{\prime }-p}d_{\omega ,p}(Y,Z)^{p},
\end{eqnarray*}%
which gives inequality (\ref{interpolation}). $\sup_{n}\left\Vert
Y(n)\right\Vert _{\omega ,p}<\infty $ and the pointwise convergence of $Y(n)$
to $Y$ implies that $\left\Vert Y\right\Vert _{\omega ,p}<\infty $. Then, by
proposition \ref{AA}, we obtain that $Y(n)$ converges uniformly to $Y$, and
we obtain our result by applying inequality (\ref{interpolation}).
\end{proof}

A similar proof gives the following proposition:

\begin{proposition}
Let $1\leq p<p^{\prime }<\infty $. Then for all $G^{m}(V)$-valued paths $%
Y,Z, $%
\begin{equation}
d_{p^{\prime }-var}(Y,Z)\leq d_{\infty }(Y,Z)^{1-p/p^{\prime
}}d_{p-var}(Y,Z)^{\frac{p}{p^{\prime }}}.  \label{interpolation2}
\end{equation}%
In particular, if $Y(n)$ converges uniformly to $Y$ and $\sup_{n}\left\Vert
Y(n)\right\Vert _{p-var}<\infty $ then 
\begin{equation*}
d_{p-var}(Y(n),Y)\rightarrow 0.
\end{equation*}
\end{proposition}

\section{The Spaces $C^{p-var}(G^{m}(V))$ and $C^{\protect\omega %
,p}(G^{m}(V))$}

We are going to prove that elements in $C^{p-var}(G^{m}(V))$ and $C^{\omega
,p}(G^{m}(V))$ are still limit in some sense of signature of smooth paths.
We will use the following proposition.

\begin{proposition}
\label{PropForMainThm}For every $g\in G^{m}(V)$, there exists a smooth path $%
h_{g}\left( .\right) =y(.)$ of constant speed $\left\vert \dot{y}\right\vert 
$ such that $S_{m}(y)_{0,1}=g$ and such that%
\begin{equation*}
\left\Vert \dot{y}\right\Vert _{L^{\infty }}\leq 2\left\Vert g\right\Vert
_{m}.
\end{equation*}%
As a consequence, 
\begin{equation*}
\left\Vert S_{m}(y)_{s,t}\right\Vert _{m}\leq 2\left\Vert g\right\Vert
_{m}(t-s).
\end{equation*}
\end{proposition}

\begin{proof}
Without lost of generalities, we can assume that $\left\Vert g\right\Vert
>0. $ (Indeed, if $\left\Vert g\right\Vert =0,$ the path $y:[0,1]\rightarrow
V,$ $u\rightarrow 0$ will do.)

By definition of the Carnot-Caratheodory norm, for every $\varepsilon >0$,
there exists a smooth path $y$, which we may take of constant speed (by time
reparametrization), such that $S_{m}(y)_{0,1}=g$ and such that the sup norm
of $\dot{y}$ is bounded by $\left\Vert g\right\Vert +\varepsilon $. Taking $%
\varepsilon =\left\Vert g\right\Vert $ finishes the first part. The second
part follows from Proposition \ref{obvProp}.
\end{proof}

\begin{remark}
We would have liked to define $h_{g}$ as a geodesic associated to $g$, e.g.
the shortest connection of the neutral element in $G^{m}(V)$ with $g$ w.r.t.
the Carnot-Caratheodory distance $d_{G^{m}(V)}$. Unfortunately, smoothness
of such geodesics is still an open problem for $m\geq 3$ (personal
communication, R.Montgomery). An affirmative answer for the case $m=2$ is
found in \cite{LS}.
\end{remark}

\begin{theorem}
\label{withH}Let $\omega $ be a control satisfying condition $(H_{p})$. A
path $Y$ belongs to $C^{\omega ,p}(G^{m}(V))$ if and only if there exists a
sequence of smooth $V$-valued paths $y(n)$ such that\newline
(i): $\sup_{n}\left\Vert S_{m}(y(n))\right\Vert _{\omega ,p}<\infty .$%
\newline
(ii): $S_{m}(y(n))$ converges pointwise to $Y$.\newline
In particular, $S_{m}(y(n))$ converges\ to $Y$ in the topology induced by $%
d_{\omega ,q}$, whenever $q>p$.
\end{theorem}

\begin{proof}
The fact that the existence of a sequence of smooth paths $y_{n}$ satisfying
conditions (i) and (ii) implies that $\left\Vert Y\right\Vert _{\omega
,p}<\infty $ is obvious. We prove the reverse implication.\newline
We let $\phi $ be a non-decreasing function in $C^{\infty }([0,1],\mathbb{R}%
) $ such that 
\begin{eqnarray*}
\phi (0) &=&0, \\
\phi (1) &=&1, \\
\forall k &\geq &1,\phi ^{(k)}(0)=\phi ^{(k)}(1)=0.
\end{eqnarray*}%
Fix a subdivision of $[0,1]$,$D=\left\{ 0=t_{0}<t_{1}<...<t_{n}=1\right\} $.
>From this subdivision, we construct a smooth path $y(D),$ one time-interval
after the other: first $y(D)_{0}=0$. Then, for $t\in \left[ t_{i-1},t_{i}%
\right] ,$ we let $y(D)_{t_{i},t}=h_{Y_{t_{i},t_{i+1}}}\left( \phi \left( 
\frac{t-t_{i}}{t_{i+1}-t_{i}}\right) \right) $ ($h_{g}$ has been defined in
the previous proposition). Thanks to our choice of the function $\phi $, $%
y(D)$ is a smooth path.

Using Proposition \ref{PropForMainThm}, for $t_{i-1}\leq s\leq t\leq t_{i}$,

\begin{eqnarray*}
\left\Vert S_{m}\left( y(D)\right) _{s,t}\right\Vert _{m} &=&\left\Vert
S_{m}\left( h_{Y_{t_{i},t_{i+1}}}\right) _{\phi \left( \frac{s-t_{i}}{%
t_{i+1}-t_{i}}\right) ,\phi \left( \frac{t-t_{i}}{t_{i+1}-t_{i}}\right)
}\right\Vert _{m} \\
&\leq &2\left\Vert Y_{t_{i},t_{i+1}}\right\Vert _{m}\left( \phi \left( \frac{%
t-t_{i}}{t_{i+1}-t_{i}}\right) -\phi \left( \frac{s-t_{i}}{t_{i+1}-t_{i}}%
\right) \right) \\
&\leq &2\left\vert \phi ^{\prime }\right\vert _{\infty }\frac{t-s}{%
t_{i+1}-t_{i}}\left\Vert Y_{t_{i},t_{i+1}}\right\Vert _{m} \\
&\leq &C\frac{t-s}{t_{i+1}-t_{i}}\left\Vert Y\right\Vert _{\omega ,p}\left[
\omega (t_{i},t_{i+1})\right] ^{1/p} \\
&\leq &C\left\Vert Y\right\Vert _{\omega ,p}\omega (s,t)^{1/p}\text{ \ \ \ \
\ \ using condition (H).}
\end{eqnarray*}%
Note that%
\begin{equation}
S_{m}\left( y(D)\right) _{t_{i}}=Y_{t_{i}}\text{ for all }i=0,...,n.
\label{ZisY}
\end{equation}%
\newline
Then for all $t_{i-1}\leq s\leq t_{i}\leq t_{j}<t\leq t_{j+1}$,%
\begin{eqnarray*}
\left\Vert S_{m}(y(D))_{s,t}\right\Vert _{m}^{p} &=&\left\Vert
S_{m}(y(D))_{s,t_{i}}\otimes Y_{t_{i},t_{j}}\otimes
S_{m}(y(D))_{t_{j},t}\right\Vert _{m}^{p} \\
&\leq &3^{p-1}\left( \left\Vert S_{m}(y(D))_{s,t_{i}}\right\Vert
_{m}^{p}+\left\Vert Y_{t_{i},t_{j}}\right\Vert _{m}^{p}+\left\Vert
S_{m}(y(D))_{t_{j},t}\right\Vert _{m}^{p}\right) \\
&\leq &C\left\Vert Y\right\Vert _{\omega ,p}^{p}\left( \omega
(s,t_{i})+\omega (t_{i},t_{j})+\omega (t_{j},t)\right) \\
&\leq &C\left\Vert Y\right\Vert _{\omega ,p}^{p}\omega (s,t).
\end{eqnarray*}%
Moreover, for $t$ as above, from equality (\ref{ZisY}) and the left
invariance of the distance $d$, we get that if $t_{j}\leq t\leq t_{j+1}$,%
\begin{eqnarray}
d_{m}(S_{m}(y(D))_{t},Y_{t}) &=&d_{m}(S_{m}(y(D))_{t_{j},t},Y_{t_{j},t}) 
\notag \\
&\leq &C\left\Vert S_{m}(y(D))_{t_{j},t}\right\Vert _{m}+\left\Vert
Y_{t_{j},t}\right\Vert _{m}  \notag \\
&\leq &C\left( \left\Vert Y_{t_{j},t_{j+1}}\right\Vert _{m}+\left\Vert
Y_{t_{j},t}\right\Vert _{m}\right)  \notag \\
&\leq &C\sup_{\substack{ s,t\in \lbrack 0,1]  \\ \left\vert t-s\right\vert
\leq \text{mesh}\left( D\right) }}\left\Vert Y_{s,t}\right\Vert _{m}.
\label{UC}
\end{eqnarray}%
$Y$ is continuous and defined on a compact ($[0,1]$), hence by
Heine-Cantor's theorem, it is uniformly continuous. Therefore, for all $%
\varepsilon >0$, there exists $\eta $ such that $\left\vert D\right\vert
<\eta \Rightarrow d_{m}(S_{m}(y(D))_{t},Y_{t})<\varepsilon $. We have just
shown that if $\left( D_{n}\right) _{n}$ is a family of subdivision of $%
[0,1] $, whose mesh goes to $0$ when $n\rightarrow \infty $, then $y(D_{n})$
is a sequence of smooth path satisfying conditions (i) and (ii).\newline
The last statement is just a corollary of inequality (\ref{interpolation}).
\end{proof}

\begin{corollary}
\label{Cor_Holder}Let $Y\mathbf{\ }$be a $G^{m}(V)$-valued path. Then $Y$ is 
$1/p$-H\"{o}lder if and only if there exists a sequence of infinitely
differentiable $V$-valued paths $y(n)$ such that\newline
(i): The $1/p$-H\"{o}lder norm of $S_{m}(y(n))$ is uniformly bounded.\newline
(ii): $S_{m}(y(n))$ converges pointwise to $Y$.\newline
In particular, given an $\alpha =1/p$ H\"{o}lder regular $G^{m}(V)$-valued
path $Y$ there is a sequence of signature of smooth paths that converge in $%
\alpha ^{\prime }$-H\"{o}lder topology to $Y$, for any $\alpha ^{\prime }$%
\thinspace $<\alpha .$
\end{corollary}

\begin{proof}
Apply the previous theorem with the control $\left( s,t\right) \mapsto t-s$.
\end{proof}

We now consider the $p$-variation distance. First, if $Y$ are paths of
finite $p$-variation, we define 
\begin{equation}
\delta _{Y}^{p}(s,t)=\sup_{D=(s\leq t_{0}<\cdots <t_{n}\leq
t)}\sum_{i=0}^{n-1}\left\Vert Y_{t_{i},t_{i+1}}\right\Vert _{m}^{p}.
\label{smallestcontrol}
\end{equation}%
In other words, $\delta _{Y}^{p}$ is the smallest control of the $p$%
-variation of $Y$.

\begin{theorem}
\label{withpvar}$Y$ belongs to $C^{p-var}(G^{m}(V))$ if and only if there
exists a sequence of infinitely differentiable $V$-valued paths $y(n)$ such
that\newline
(i): $\left\Vert S_{m}(y(n))\right\Vert _{p-var}$ is uniformly bounded.%
\newline
(ii): $S_{m}(y(n))$ converges pointwise to $Y$.\newline
In particular, $S_{m}(y(n))$ converges in $q$-variation to $Y$, whenever $%
q>p $.
\end{theorem}

\begin{proof}
We construct $y(D)$ from $Y$ as in the proof of theorem \ref{withH} with the
help of a subdivision $D=\left\{ 0=t_{0}<t_{1}<...<t_{n}=1\right\} $. Define
the control%
\begin{equation*}
\omega _{D}(s,t)=\left( \frac{t-s}{t_{i+1}-t_{i}}\right) ^{p}\delta
_{Y}^{p}(t_{i},t_{i+1})\text{ for }t_{i}\leq s\leq t\leq t_{i+1},\text{ }%
0\leq i\leq n-1;
\end{equation*}%
and for $0\leq i<j\leq n-1$, and $t_{i-1}\leq s\leq t_{i}\leq t_{j}\leq
t\leq t_{j+1},$ 
\begin{equation}
\omega _{D}(s,t)=\omega _{D}(s,t_{i})+\delta _{Y}^{p}(t_{i},t_{j})+\omega
_{D}(t_{j},t).  \label{omegaD}
\end{equation}%
It is easy to check that $\omega _{D}$ is a control (but does not
necessarily satisfies condition $(H_{p})$). Then, from the proof of theorem %
\ref{withH}, we see that for $t_{i}\leq s\leq t\leq t_{i+1}$,%
\begin{eqnarray*}
\left\Vert S_{m}(y(D))_{s,t}\right\Vert ^{p} &\leq &C\left( \frac{t-s}{%
t_{i+1}-t_{i}}\right) ^{p}\delta _{Y}^{p}(t_{i},t_{i+1}) \\
&\leq &C\omega _{D}(s,t).
\end{eqnarray*}%
Then,\ if $t_{i-1}\leq s\leq t_{i}\leq t_{j}<t\leq t_{j+1}$,%
\begin{eqnarray*}
\left\Vert S_{m}(y(D))_{s,t}\right\Vert _{m}^{p} &\leq &C\left( \left\Vert
S_{m}(y(D))_{s,t_{i}}\right\Vert _{m}^{p}+\left\Vert
Y_{t_{i},t_{j}}\right\Vert _{m}^{p}+\left\Vert
S_{m}(y(D))_{t_{j},t}\right\Vert _{m}^{p}\right) \\
&\leq &C\left( \omega _{D}(s,t_{i})+\delta _{Y}^{p}(t_{i},t_{j})+\omega
_{D}(t_{j},t)\right) \\
&=&C\omega _{D}(s,t).
\end{eqnarray*}%
Therefore,%
\begin{equation*}
\left\Vert S_{m}(y(D))\right\Vert _{p-var}\leq C\omega
_{D}(0,1)^{1/p}=C\delta _{Y}^{p}(0,1)^{1/p}=C\left\Vert Y\right\Vert
_{p-var}.
\end{equation*}%
\newline
Hence, if $\left( D_{n}=\left\{ 0\leq t_{1}^{n}<\cdots \leq
t_{\#D_{n}}^{n}\right\} \right) _{n}$ is a sequence of subdivision of $[0,1]$
whose mesh tends to $0$, we have just proved that $\left( y(D_{n})\right)
_{n}$ satisfies condition (i); condition (ii) is treated just like before,
using inequality (\ref{UC}).\newline
The last statement is just a corollary of inequality \ref{interpolation2},
once we prove that $S(y(D_{n}))$ converges uniformly to $Y.$ To do so,
define $h_{D_{n}}(\delta )=\sup_{\left\vert t-s\right\vert \leq \delta
}\omega _{D_{n}}(s,t)$, and $h_{\infty }(\delta )=\sup_{\left\vert
t-s\right\vert \leq \delta }\omega (s,t)$. By Heine-Cantor's theorem, $%
\omega $ is uniformly continuous. As it is zero on the diagonal, we obtain
that $h_{D_{n}}(\delta )\rightarrow _{\delta \rightarrow 0}0$. If $%
\left\vert D_{n}\right\vert \leq \delta $, for all $s<t$ such that $%
\left\vert t-s\right\vert \leq \delta $, there exists $1\leq i\leq \#D_{n}-2$
such that $t_{i}\leq s<t\leq t_{i+2}$. Hence, by definition of 
\begin{equation*}
\omega _{D_{n}}(s,t)\leq \omega (t_{i}^{n},t_{i+2}^{n})\leq h_{\infty
}(2\left\vert D_{n}\right\vert )\leq h_{\infty }(2\delta ).
\end{equation*}%
Hence, given an $\varepsilon >0$, there exists $\delta _{0}>0$ such that $%
\delta \leq \delta _{0}\Rightarrow h_{\infty }(2\delta )\leq \varepsilon .$
Let $N$ be such that $\sup_{n\geq N}\left\vert D_{n}\right\vert \leq \delta
_{0}$. Then, there exists $\delta _{1}>0$ such that $\delta \leq \delta
_{1}\Rightarrow \max_{n<N}h_{D_{n}}(\delta )\leq \varepsilon $. In
particular, $\delta \leq \min \left\{ \delta _{0},\delta _{1}\right\}
\Rightarrow \sup_{n\in \mathbb{N}\cup \left\{ \infty \right\}
}h_{D_{n}}(\delta )\leq \varepsilon $, i.e. $\sup_{n\in \mathbb{N}\cup
\left\{ \infty \right\} }h_{D_{n}}$ goes to $0$ at $0$. We have therefore
proved that $S(y(D_{n}))$ is equicontinuous, which implies the wanted
uniform convergence by lemma \ref{AA}.
\end{proof}

\begin{remark}
We could have obtained the above theorem replacing smooth paths by paths of
bounded variation, by a simple change of time. Indeed, a path of finite $p$%
-variation can be reparametrized into a $1/p$-H\"{o}lder path.
\end{remark}

\begin{remark}
These convergences results, more precisely the (almost) geodesic
approximation versus piecewise linear, can be compared with probabilistic
constructions in the context of Brownian rough paths (\cite{LQ,F,FV}). There
it is essential to base approximations on nested (or dyadic) subdivisions of 
$[0,1]$ which is not required here.
\end{remark}

\section{The Spaces $C^{0,p-var}(G^{m}(V))$ and $C^{0,\protect\omega %
,p}(G^{m}(V))$}

Recall the definitions of these space made earlier. We are first going to
give an equivalent definition of these sets. We will then prove that these
spaces are separable.

\subsection{A Ciesielski/Museliak-Semadini Type Result}

Ciesielski and Museliak-Semadini proved at similar times with different
techniques the following theorem, in the case of H\"{o}lder real valued
paths. Taking $\omega (s,t)=t-s$, and $m=1,$ $V=\mathbb{R}$ in the theorem
below gives (some) of their results. \cite{Ci,MS}.

\begin{theorem}
Let $Y$ be an element of $C^{\omega ,p}(G^{m}(V))$. We assume that $\omega $
satisfies condition $(H_{p})$ and that 
\begin{equation}
\lim_{\delta \rightarrow 0}\sup_{\substack{ 0\leq s<t\leq 1  \\ t-s\geq
\delta }}\frac{t-s}{\omega (s,t)^{1/p}}=0.  \label{condHPP}
\end{equation}%
Then, $Y$ belongs to $C^{0,\omega ,p}(G^{m}(V))$ if and only if%
\begin{equation*}
\lim_{\delta \rightarrow 0}\sup_{\substack{ 0\leq s<t\leq 1  \\ t-s\leq
\delta }}\frac{\left\Vert Y_{s,t}\right\Vert }{\omega (s,t)^{1/p}}=0\text{.}
\end{equation*}
\end{theorem}

Note that the condition (\ref{condHPP}) implies that $p>1$. Under condition $%
\left( H_{p}\right) $, we have%
\begin{equation*}
\overline{\lim_{\delta \rightarrow 0}}\sup_{\substack{ 0\leq s<t\leq 1  \\ %
t-s\geq \delta }}\frac{t-s}{\omega (s,t)^{1/p}}=:C<\infty .
\end{equation*}%
Then condition \ref{condHPP} simply reads $C=0.$

\begin{proof}
Assume that $Y$ belongs to $C^{0,\omega ,p}(G^{m}(V))$. Then, by definition,
there exists a sequence of signature of smooth paths $(S_{m}(y_{n}))_{n}$
such that 
\begin{equation*}
\lim_{n\rightarrow \infty }d_{\omega ,p}(S_{m}(y_{n}),Y)=0.
\end{equation*}%
For all $s<t$, 
\begin{equation*}
\frac{\left\Vert Y_{s,t}\right\Vert }{\omega (s,t)^{1/p}}\leq d_{p,\omega
}(Y,S_{m}(y_{n}))+\frac{\left\Vert S_{m}(y_{n})_{s,t}\right\Vert }{\omega
(s,t)^{1/p}},
\end{equation*}%
hence for all $n$, 
\begin{equation*}
\lim_{\delta \rightarrow 0}\sup_{\substack{ 0\leq s<t\leq 1  \\ t-s\geq
\delta }}\frac{\left\Vert Y_{s,t}\right\Vert }{\omega (s,t)^{1/p}}\leq
d_{p,\omega }(Y,S_{m}(y_{n}))+\lim_{\delta \rightarrow 0}\sup_{\substack{ %
0\leq s<t\leq 1  \\ t-s\leq \delta }}\frac{\left\Vert
S_{m}(y_{n})_{s,t}\right\Vert }{\omega (s,t)^{1/p}}.
\end{equation*}%
But as $y_{n}$ is smooth, $S_{m}(y_{n})$ is Lipschitz, hence, by the
assumption on the control $\omega $, we obtain that $\lim_{\delta
\rightarrow 0}\sup_{\substack{ 0\leq s<t\leq 1  \\ t-s\leq \delta }}\frac{%
\left\Vert S_{m}(y_{n})_{s,t}\right\Vert }{\omega (s,t)^{1/p}}=0$.
Therefore, for all $n$, 
\begin{equation*}
\lim_{\delta \rightarrow 0}\sup_{\substack{ 0\leq s<t\leq 1  \\ t-s\leq
\delta }}\frac{\left\Vert Y_{s,t}\right\Vert }{\omega (s,t)^{1/p}}\leq
d_{p,\omega }(Y,S_{m}(y_{n})),
\end{equation*}%
i.e. this limit is equal to $0$.

Reciprocally, assume that 
\begin{equation*}
\lim_{\delta \rightarrow 0}\sup_{\substack{ 0\leq s<t\leq 1  \\ t-s\leq
\delta }}\frac{\left\Vert Y_{s,t}\right\Vert }{\omega (s,t)^{1/p}}=0\text{.}
\end{equation*}%
We define%
\begin{equation*}
\vartheta _{Y}(\delta )=\sup_{\substack{ 0\leq s<t\leq 1  \\ t-s\leq \delta 
}}\frac{\left\Vert Y_{s,t}\right\Vert }{\omega (s,t)^{1/p}}\in \lbrack
0,\left\Vert Y\right\Vert _{\omega ,p}]\text{ for }\delta \in \lbrack 0,1].
\end{equation*}%
We now define $y(D)$ from $Y$ as in the proof of theorem \ref{withH}, where $%
D$ is a given subdivision of $[0,1]$. With techniques similar to the one
used in the proof of theorem \ref{withH}, we see that $\vartheta
_{S_{m}(y(D))}(\delta )\leq C\left\Vert Y\right\Vert _{\omega ,p}\vartheta
_{Y}(\delta )$ for a universal constant $C$. Then, for $s<t$ such that $%
\left\vert t-s\right\vert \geq $ $\delta $,%
\begin{eqnarray*}
\frac{d(Y_{s,t},S_{m}(y(D))_{s,t})}{\omega (s,t)^{1/p}} &\leq &\frac{%
d_{\infty }(Y,S_{m}(y(D)))}{\inf_{\left\vert t-s\right\vert \geq \delta
}\omega (s,t)^{1/p}} \\
&\leq &C\frac{d_{\infty }(Y,S_{m}(y(D)))}{\delta },\text{ \ \ \ using
condition }\left( H_{p}\right)
\end{eqnarray*}%
For $s<t$ such that $\left\vert t-s\right\vert <\delta ,$%
\begin{eqnarray*}
\frac{d(Y_{s,t},S_{m}(y(D))_{s,t})}{\omega (s,t)^{1/p}} &\leq &\frac{%
\left\Vert Y_{s,t}\right\Vert _{m}+\left\Vert S_{m}(y(D))_{s,t}\right\Vert
_{m}}{\omega (s,t)^{1/p}} \\
&\leq &C\left\Vert Y\right\Vert _{\omega ,p}\vartheta (\delta ).
\end{eqnarray*}%
Hence, for all $\delta >0$,%
\begin{equation*}
d_{\omega ,p}(Y,S_{m}(y(D)))\leq C\max \left\{ \frac{d_{\infty
}(Y,S_{m}(y(D)))}{\delta },\left\Vert Y\right\Vert _{\omega ,p}\vartheta
(\delta )\right\} .
\end{equation*}%
Now, once again, consider a sequence of subdivisions $\left( D_{n}\right)
_{n}$ of $[0,1]$ whose mesh tends to $0$. For a given $\varepsilon >0$,
there exists $\delta >0$ such that $\left\Vert Y\right\Vert _{\omega
,p}\vartheta (\delta )\leq \varepsilon /C$ (as by assumption, $\lim_{\delta
\rightarrow 0}\vartheta (\delta )=0$). For such a $\delta $ and $\varepsilon 
$, there exists $N$ such that if $n\geq N$, $\frac{d_{\infty
}(Y,S_{m}(y(D_{n})))}{\delta }\leq \varepsilon /C$. We have therefore proved
that for all $\varepsilon >0$, there exists $N$ such that for all $n\geq N,$ 
$d_{\omega ,p}(Y,S_{m}(y(D_{n})))\leq \varepsilon $.
\end{proof}

\subsection{A Wiener Type Result}

We now prove a similar theorem, but in $p$-variation topology rather than
modulus topology. $p$-variation closure of step functions has been
characterized by Wiener \cite{Wi,MS}, for real valued functions, but not
necessarily continuous. We obtain this result here in the simpler case of
continuous paths, but harder case of group valued paths.

\begin{theorem}
\label{Wiener}Let $Y$ be an element of $C^{p-var}(G^{m}(V))$, $p>1$. Then, $%
Y $ belongs to $C^{0,p-var}(G^{m}(V))$ if and only if%
\begin{equation*}
\lim_{\delta \rightarrow 0}\sup_{\substack{ D=(0=t_{0}<\cdots <t_{n}=1)  \\ %
\left\vert D\right\vert \leq \delta }}\sum_{i=0}^{n-1}\delta
_{Y}^{p}(t_{i},t_{i+1})=0,
\end{equation*}%
with $\delta _{Y}^{p}$ defined as earlier (equation (\ref{smallestcontrol})).
\end{theorem}

\begin{proof}
If $Y$ belongs to $C^{0,p-var}(G^{m}(V))$, there exists a sequence of smooth
paths $y_{n}$ such that $S_{m}(y_{n})$ converges in the topology induced by $%
d_{p-var}$ to $Y$. Then, if $D=(0=t_{0}<\cdots <t_{l}=1)$ is a partition of $%
[0,1]$,%
\begin{equation*}
\left( \sum_{i=0}^{l-1}\delta _{Y}^{p}(t_{i},t_{i+1})\right) ^{1/p}\leq
\left( \sum_{i=0}^{l-1}\delta _{S_{m}(y_{n})}^{p}(t_{i},t_{i+1})\right)
^{1/p}+d_{p-var}(Y,S_{m}(y_{n})).
\end{equation*}%
Because $y_{n}$ is smooth and $p>1$, 
\begin{equation*}
\lim_{\delta \rightarrow 0}\sup_{\substack{ D=(0=t_{0}<\cdots <t_{l}=1)  \\ %
\left\vert D\right\vert \leq \delta }}\sum_{i=0}^{l-1}\omega
_{p,S_{m}(y_{n})}(t_{i},t_{i+1})=0.
\end{equation*}%
We therefore obtain that 
\begin{equation*}
\lim_{\delta \rightarrow 0}\sup_{\substack{ D=(0=t_{0}<\cdots <t_{l}=1)  \\ %
\left\vert D\right\vert \leq \delta }}\sum_{i=0}^{l-1}\delta
_{Y}^{p}(t_{i},t_{i+1})=0.
\end{equation*}

Reciprocally, we let $\pi _{n}=\left\{ k2^{-n},k\in \left\{ 0,\cdots
,2^{n}\right\} \right\} $ be a (very specific for simplicity) sequence of
subdivisions of $[0,1]$, whose mesh goes to $0$. We let $y_{n}=y(\pi _{n})$,
obtained from $\pi _{n}$ and $Y,$ as in the proof of theorems \ref{withpvar}%
. We also define $\omega _{n}(s,t)=\omega _{\pi _{n}}(s,t)$ ($\omega _{\pi
_{n}}$ is defined as in equation (\ref{omegaD})), $\omega _{\infty }=\delta
_{Y}^{p}$ and for $n\in \mathbb{N\cup }\left\{ \mathbb{\infty }\right\} $,%
\begin{equation*}
g_{n}(\delta )=\sup_{\substack{ D=(0=t_{0}<\cdots <t_{l}=1)  \\ \left\vert
D\right\vert \leq \delta }}\sum_{i=0}^{l-1}\omega _{n}(t_{i},t_{i+1}).
\end{equation*}%
By assumption, $g_{\infty }$ tends to $0$ at $0$. We are now going to prove
that $g_{n}$, $n\in \mathbb{N}$, share the same property. Consider a
subdivision $D=(0=t_{0}<\cdots <t_{l}=1)$ of $[0,1]$ with mesh size less
than $\delta <2^{-n}$. From the definition of $\omega _{n}$, we see that
when we compute $\sum_{i=0}^{l-1}\omega _{n}(t_{i},t_{i+1})$, we obtain the
same result if we add to our subdivision $D$ all the point $\left\{
k2^{-n},0\leq k\leq 2^{n}\right\} .$ Having done so, we can write%
\begin{eqnarray}
\sum_{i=0}^{l-1}\omega _{n}(t_{i},t_{i+1}) &=&\sum_{k=0}^{2^{n}-1}\omega
_{\infty }\left( \frac{k}{2^{n}},\frac{k+1}{2^{n}}\right) \sum_{i,\text{ }%
\frac{k}{2^{n}}\leq t_{i}\leq t_{i+1}\leq \frac{k+1}{2^{n}}}\left( \frac{%
t_{i+1}-t_{i}}{2^{-n}}\right) ^{p}  \notag \\
&\leq &\sum_{k=0}^{2^{n}-1}\omega _{\infty }\left( \frac{k}{2^{n}},\frac{k+1%
}{2^{n}}\right) \sum_{i,\text{ }\frac{k}{2^{n}}\leq t_{i}\leq t_{i+1}\leq 
\frac{k+1}{2^{n}}}\left( \frac{t_{i+1}-t_{i}}{2^{-n}}\right) \left( \frac{%
\delta }{2^{-n}}\right) ^{p-1}  \notag \\
&=&\left( \frac{\delta }{2^{-n}}\right) ^{p-1}\sum_{k=0}^{2^{n}-1}\omega
_{\infty }\left( \frac{k}{2^{n}},\frac{k+1}{2^{n}}\right)  \notag \\
&\leq &\left( \frac{\delta }{2^{-n}}\right) ^{p-1}g_{\infty }(2^{-n}).
\label{subdivlessthan2_n}
\end{eqnarray}%
That proves that $\lim_{\delta \rightarrow 0}g_{n}(\delta )=0$. We actually
claim that $\lim_{\delta \rightarrow 0}\sup_{n\in \mathbb{N\cup }\left\{ 
\mathbb{\infty }\right\} }g_{n}(\delta )=0$.\newline
Again, fix a subdivision $D=(0=t_{0}<\cdots <t_{l}=1)$ of $[0,1]$ . Assume
first that the mesh of $D$ is greater than or equal to $2^{-n}$. Let $%
I=\left\{ 0\leq i\leq 2^{n},D\cap \lbrack i2^{-n},(i+1)2^{-n})\neq
\varnothing \right\} $, and $s_{i}$ the smallest element of $D\cap \lbrack
i2^{-n},(i+1)2^{-n})$ for $i\in I.$ We also let $I=\left\{ i_{1},\cdots
,i_{\left\vert I\right\vert }\right\} $. $D^{\prime }=\left(
0=s_{i_{1}}<\cdots <s_{i_{\left\vert I\right\vert }}=1\right) $ is then a
subdivision of $[0,1]$, included in $D$, such that for all $1\leq j\leq
\left\vert I\right\vert $, $\left\vert s_{i_{j+1}}-s_{i_{j}}\right\vert \leq
\left\vert D\right\vert +2^{-n}\leq 2\left\vert D\right\vert $. In
particular, this implies, as $\frac{i_{j}}{2^{n}}\leq s_{i_{j}}\leq \frac{%
i_{j}+1}{2^{n}},$ that 
\begin{equation*}
i_{j+1}2^{-n}-\left( i_{j}+1\right) 2^{-n}\leq s_{i_{j+1}}-s_{i_{j}}\leq
2\left\vert D\right\vert .
\end{equation*}%
By the super-additivity property of a control, $\sum_{i=0}^{l-1}\omega
_{n}(t_{i},t_{i+1})\leq \sum_{j=1}^{\left\vert I\right\vert -1}\omega
_{n}(s_{i_{j}},s_{i_{j+1}})$. Moreover,%
\begin{eqnarray*}
\sum_{j=1}^{\left\vert I\right\vert -1}\omega _{n}(s_{i_{j}},s_{i_{j+1}})
&=&\sum_{j=1}^{\left\vert I\right\vert -1}\left( \omega
_{n}(s_{i_{j}},(i_{j}+1)2^{-n})+\omega _{n}(i_{j+1}2^{-n},s_{i_{j+1}})\right)
\\
&&+\sum_{j=1}^{\left\vert I\right\vert -1}\omega _{\infty }(\left(
i_{j}+1\right) 2^{-n},i_{j+1}2^{-n})\text{ \ \ \ by definition of }\omega
_{n} \\
&=&\sum_{j=1}^{\left\vert I\right\vert -1}\omega
_{n}(s_{i_{j}},(i_{j}+1)2^{-n})+\sum_{j=2}^{\left\vert I\right\vert }\omega
_{n}(i_{j}2^{-n},s_{i_{j}}) \\
&&+\sum_{j=1}^{\left\vert I\right\vert -1}\omega _{\infty }(\left(
i_{j}+1\right) 2^{-n},i_{j+1}2^{-n}) \\
&\leq &\sum_{j=1}^{\left\vert I\right\vert -1}\left( \omega _{\infty
}(i_{j}2^{-n},(i_{j}+1)2^{-n})+\omega _{\infty }(\left( i_{j}+1\right)
2^{-n},i_{j+1}2^{-n})\right) \\
&\leq &g_{\infty }\left( 2\left\vert D\right\vert \right) .
\end{eqnarray*}%
When the mesh of $D$ is less than or equal to $2^{-n}$, we have already seen
(equation (\ref{subdivlessthan2_n})) that 
\begin{equation*}
\sum_{i=0}^{l-1}\omega _{n}(t_{i},t_{i+1})\leq \left( \frac{\left\vert
D\right\vert }{2^{-n}}\right) ^{p-1}g_{\infty }\left( 2^{-n}\right) \leq
g_{\infty }\left( 2^{-n}\right) .
\end{equation*}%
Therefore, $g_{n}(\delta )\leq \max \left\{ g_{\infty }\left( 2\delta
\right) ,g_{\infty }\left( 2^{-n}\right) \right\} $. Let $\varepsilon >0$.
There exists $n_{0}$ such that $n\geq n_{0}\Longrightarrow g_{\infty }\left(
2^{-n}\right) \leq \varepsilon $. Then, using the observation above that%
\begin{equation*}
\lim_{\delta \rightarrow 0}g_{n}(\delta )=0
\end{equation*}%
for fixed $n$, there exists $\delta _{0}$ such that for all $\delta \leq
\delta _{0}$, 
\begin{equation*}
\max \left\{ g_{\infty }\left( 2\delta \right) ,g_{n}\left( \delta \right)
,n=0,\cdots ,n_{0}-1\right\} \leq \varepsilon .
\end{equation*}%
In particular, for all $\delta \leq \delta _{0}$, $\sup_{n\in \mathbb{N\cup }%
\left\{ \mathbb{\infty }\right\} }g_{n}(\delta )\leq \varepsilon $, as%
\begin{equation*}
\sup_{n}g_{n}(\delta )=\max \left\{ \max_{n<n_{0}}g_{n}(\delta ),\sup_{n\geq
n_{0}}g_{n}(\delta )\right\} .
\end{equation*}%
In other words,%
\begin{equation*}
\lim_{\delta \rightarrow 0}\sup_{n\in \mathbb{N\cup }\left\{ \mathbb{\infty }%
\right\} }g_{n}(\delta )=0.
\end{equation*}

We consider again a subdivision $D=(0=t_{0}<\cdots <t_{l}=1)$ of $[0,1]$. As%
\begin{equation*}
d(Y_{t_{i},t_{i+1}},S_{m}(y_{n})_{t_{i},t_{i+1}})^{p}\leq \min \left\{
d_{\infty }(S_{m}(y_{n}),Y)^{p},2^{p-1}\left( \omega
_{n}(t_{i},t_{i+1})+\omega _{\infty }(t_{i},t_{i+1})\right) \right\} ,
\end{equation*}%
we obtain that 
\begin{eqnarray*}
\sum_{i=0}^{l-1}d(Y_{t_{i},t_{i+1}},S_{m}(y_{n})_{t_{i},t_{i+1}})^{p} &\leq
&\sum_{\substack{ i\in \{0,\cdots ,l-1\}  \\ \left\vert
t_{i+1}-t_{i}\right\vert >\delta }}d_{\infty }(S_{m}(y_{n}),Y)^{p} \\
&&+2^{p-1}\sum_{\substack{ i\in \{0,\cdots ,n-1\}  \\ \left\vert
t_{i+1}-t_{i}\right\vert \leq \delta }}\left( \omega
_{n}(t_{i},t_{i+1})+\omega _{\infty }(t_{i},t_{i+1})\right) \\
&\leq &\frac{d_{\infty }(S_{m}(y_{n}),Y)}{\delta }^{p}+2^{p}\sup_{n\in 
\mathbb{N\cup }\left\{ \mathbb{\infty }\right\} }g_{n}(\delta ).
\end{eqnarray*}%
Therefore, $d_{p-var}\left( Y,S_{m}(y_{n})\right) ^{p}\leq $ $\frac{%
d_{\infty }(S_{m}(y_{n}),Y)^{p}}{\delta }+2^{p}\sup_{n\in \mathbb{N\cup }%
\left\{ \mathbb{\infty }\right\} }g_{n}(\delta )$. That gives us our result.
\end{proof}

One can then see the following equality of spaces:%
\begin{eqnarray*}
C^{0,\omega ,p}(G^{m}(V)) &=&\overline{\cup _{q>p}C^{0,\omega ,q}\left(
G^{m}(V)\right) }\text{, where the closure is the }d_{p,\omega }\text{%
-closure,} \\
C^{0,p-var}(G^{m}(V)) &=&\overline{\cup _{q>p}C^{0,q-var}\left(
G^{m}(V)\right) }\text{, where the closure is the }d_{p-var}\text{-closure.}
\end{eqnarray*}

\subsection{Polishness}

\subsubsection{Separability}

\begin{theorem}
We assume that $\omega $ is such that for all $s<t$, $\omega (s,t)\geq
K\left( t-s\right) ^{p}$.\footnote{%
which is true if $\omega $ satisfies condition $\left( H_{p}\right) $.}
Then, $C^{0,\omega ,p}(G^{m}(V))$ is separable.
\end{theorem}

\begin{proof}
We know that the space $C_{0}^{1}([0,1],V)$ of continuously differentiable
paths is separable. Let $D$ be a countable set of $C_{0}^{1}([0,1],V)$ such
that its Lipschitz closure is $C_{0}^{1}([0,1],V)$. We claim that the $%
d_{\omega ,p}$-closure of $S_{m}(D)=\left\{ S_{m}\left( y\right) ,y\in
D\right\} $ is dense in $C^{0,\omega ,p}(G^{m}(V))$. Indeed, if $Y\in
C^{0,\omega ,p}(G^{m}(V))$, there exists a sequence $\left( y_{n}\right)
_{n} $ of elements in $C_{0}^{1}([0,1],V)$, such that 
\begin{equation*}
d_{\omega ,p}\left( S_{m}\left( y_{n}\right) ,Y\right) \rightarrow
_{n\rightarrow \infty }0.
\end{equation*}%
For all $n$, there exists $\widetilde{y_{n}}\in D$ such that the Lipschitz
distance between $y_{n}$ and $\widetilde{y_{n}}$ goes to $0$ with $%
n\rightarrow \infty .$ By theorem 1 with its continuity statement in \cite%
{Ly}, we deduce that the Lipschitz distance between $S\left( y_{n}\right) $
and $S\left( \widetilde{y_{n}}\right) $ goes to $0$ when $n$ tends to
infinity. In particular, $d_{\omega ,p}\left( S\left( y_{n}\right) ,S\left( 
\widetilde{y_{n}}\right) \right) $ tends to $0$ when $n\rightarrow \infty $.
The triangle inequalities show that $S\left( \widetilde{y_{n}}\right) $
converges to $Y$ in the topology induced by $d_{\omega ,p}$.
\end{proof}

The same proof gives the following theorem:

\begin{theorem}
$C^{0,p-var}(G^{m}(V))$ is separable.
\end{theorem}

We now look at the space $C^{p-var}(G^{m}(V))$ and $C^{\omega ,p}(G^{m}(V)$.

\begin{theorem}
$C^{p-var}(G^{m}(V))$ and $C^{\omega ,p}(G^{m}(V))$ are not separable.
\end{theorem}

\begin{proof}
The proof is pretty simple. If they were separable, it would mean,
projecting $G^{m}(V)$ onto $V$, that $C^{p-var}(V)$ and $C^{\omega ,p}(V)$
are separable. They are not, see \cite{Ci,MS}.
\end{proof}

\subsubsection{Completeness}

\begin{theorem}
$C^{p-var}(G^{m}(V))$ is complete.
\end{theorem}

\begin{proof}
We suppose that $X^{n}$ is a Cauchy sequence. Then, it is also a Cauchy
sequence for the sup norm, therefore, it converges when $n$ tends to
infinity, in sup norm to, say, $X$. For a given $\varepsilon >0$, there
exists $N$, such that $n,m\geq N$ implies that for all subdivision $D=\left(
0\leq t_{0}<t_{1}<\cdots <t_{l}\leq 1\right) $ of $[0,1]$,%
\begin{equation*}
\sum_{i=0}^{l-1}d\left( X_{t_{i},t_{i+1}}^{n},X_{t_{i},t_{i+1}}^{m}\right)
^{p}<\varepsilon .
\end{equation*}%
In particular, letting $m$ tends to infinity, we obtain that $%
\sum_{i=0}^{l-1}d\left( X_{t_{i},t_{i+1}}^{n},X_{t_{i},t_{i+1}}\right)
^{p}<\varepsilon $, this being true for all subdivisions. That proves our
assertion.
\end{proof}

A similar and somewhat simpler proof gives the following:

\begin{theorem}
Let $\omega $ be a control which is zero only on the diagonal. Then $%
C^{\omega ,p}(G^{m}(V))$ is complete.
\end{theorem}

Of course, $C^{0,\omega ,p}(G^{m}(V))$ and $C^{0,p-var}(G^{m}(V))$ are
complete, being closed subsets of complete sets.

\section{Conclusion}

We have therefore characterized precisely the spaces $C^{0,p-var}(G^{m}(V))$%
, $C^{p-var}(G^{m}(V))$, $C^{0,\omega ,p}(G^{m}(V))$ and $C^{0,\omega
,p}(G^{m}(V))$. The separability statement, incidentally, proves that the
inclusions of $C^{0,p-var}(G^{m}(V))$ in $C^{p-var}(G^{m}(V))$ and $%
C^{0,p-var}(G^{m}(V))$ in $C^{p-var}(G^{m}(V))$ are\ strict. We let $\omega
(s,t)=t-s$ for the rest of this discussion. The function $g$%
\begin{eqnarray*}
\mathbb{R} &\rightarrow &\mathbb{R} \\
t &\rightarrow &\sum_{i=1}^{\infty }2^{-i/p}\sin (2^{i}t)
\end{eqnarray*}%
provides a concrete proof of the strict inclusion of $C^{0,p-var}(G^{m}(V))$
in $C^{p-var}(G^{m}(V))$ and $C^{0,p-var}(G^{m}(V))$ in $C^{p-var}(G^{m}(V))$
\cite{Du}. Also, if $V\mathbb{=R}^{2}$ is generated by a basis $\left\{
e_{1},e_{2}\right\} $, $X_{t}=\exp \left( t[e_{1},e_{2}]\right) $ is an
element of $C^{\omega ,2}(G^{m}(V))\backslash C^{0,\omega ,2}(G^{m}(V))$ and 
$C^{2-var}(G^{m}(V))\backslash C^{0,2-var}(G^{m}(V))$. The function $h$%
\begin{eqnarray*}
\mathbb{R}^{+} &\rightarrow &\mathbb{R} \\
t &\rightarrow &\left\{ 
\begin{array}{l}
\frac{t^{1/p}}{\log t}\cos ^{2}\left( \frac{\pi }{t}\right) \text{ if }t>0,
\\ 
0\text{ if }t=0.%
\end{array}%
\right.
\end{eqnarray*}%
proves that the inclusions%
\begin{eqnarray*}
\cup _{q<p}C^{q-var}(G^{m}(V)) &\subset &C^{0,p-var}(G^{m}(V)), \\
\cup _{q<p}C^{\omega ,q}(G^{m}(V)) &\subset &C^{0,\omega ,p}(G^{m}(V))
\end{eqnarray*}%
are strict \cite{Du}.

$C^{0,\omega ,p}(G^{m}(V))$ and $C^{0,p-var}(G^{m}(V))$ are therefore Polish
space. In particular, the (Stratonovich enhanced) Brownian motion takes
values in a Polish space. Many important probabilistic theorems (e.g.
Prohorov's theorem) rely on Polishness.

Finally, we want to point out the fact that the approximations of a rough
path $X$ that we have introduced (the (almost) geodesic one) may be very
useful in various area. For example, in the field of stochastic numerical
analysis, consider $\mathbf{B}$ the Stratonovich enhanced Brownian Motion
lying above a standard $d$-dimensional Brownian motion $B.$ Then let $B^{n}$
be the geodesic approximation based on the subdivision $\pi (n)=\left\{ 
\frac{k}{n},\text{ }k=0,\cdots ,n\right\} $ of $[0,1]$, and $B^{(n)}$ the
linear path which coincides with $B$ at the points $\frac{k}{n}$, and linear
in the intervals $[\frac{k-1}{n},\frac{k}{n}]$, $k=0,\cdots ,n$. Consider
the Stratonovich differential equation%
\begin{equation*}
dX_{t}=V_{0}(t,X_{t})dt+V(t,X_{t})\circ dB_{t}
\end{equation*}%
and its approximations%
\begin{eqnarray*}
dX_{t}^{n} &=&V_{0}(t,X_{t}^{n})dt+V(t,X_{t}^{n})dB_{t}^{n}, \\
dX_{t}^{\left( n\right) } &=&V_{0}(t,X_{t}^{\left( n\right)
})dt+V(t,X_{t}^{\left( n\right) })dB_{t}^{\left( n\right) }.
\end{eqnarray*}%
Then, in the $L^{2}$ sense, $X^{n}$ converges to $X$ with a speed of
convergence proportional to $\frac{1}{\sqrt{n}}$, while $X^{\left( n\right)
} $ converges to $X$ with a speed of convergence proportional to $\frac{1}{n}%
.$ Lifting $\mathbf{B}$ to a $G^{(m)}(R^{d})$-valued process and considering
the almost geodesic approximation in this group would lead to a speed of
convergence proportional to $n^{-m/2}$. See \cite{CG,KP} for a discussion
involving speed of convergence for algorithm involving the use of iterated
integrals.

\begin{acknowledgement}
The authors wish to thank R.\ Montgomery for pointing out that smoothness of
geodesics in $G^{m}(\mathbb{R}^{d}),m>2,$ is still an open problem and M.
Gubinelli for helping us to improve the presentation.
\end{acknowledgement}

\end{document}